# The Estimation of Approximation Error using the Inverse Problem and the Set of Numerical Solutions

Alekseev A.K., Bondarev A.E.
Keldysh Institute of Applied Mathematics RAS
Moscow, Russia
e-mail: aleksey.k.alekseev@gmail.com, bond@keldysh.ru

**Abstract.** The Inverse Problem for the estimation of a point-wise approximation error occurring at the discretization and solving of the system of partial differential equations is addressed. The set of the differences between the numerical solutions is used as the input data. The analyzed solutions are obtained by the numerical algorithms of the distinct inner structure on the same grid. The approximation error is estimated by the Inverse Problem, which is posed in the variational statement with the zero order Tikhonov regularization. The numerical tests, performed for the two dimensional inviscid compressible flows corresponding to Edney-I and Edney-VI shock wave interference modes, are provided. The analyzed flowfields are computed using ten different numerical algorithms. The comparison of the estimated approximation error and the true error, obtained by subtraction of numerical and analytic solutions, is presented.

**Keywords:** approximation error, differences of numerical solutions, Inverse Problem, Tikhonov regularization, Euler equations.

## 1. Introduction

The need for the reliable numerical methods at modeling of problems governed by systems of partial differential equations (PDE) causes the interest to the verification of software and numerical solutions that stimulates the development of methods for the approximation error estimation. There exists significant number of approaches for the evaluation of the approximation error [1-11]. However, some of them are highly computationally expensive; others are limited by the domain of application or provide only qualitative results. So, the computationally inexpensive, nonintrusive universal methods for the discretization error estimation are of the current interest. Herein we consider an approach to the approximation error estimation, which is based on Inverse Problem treatment of the numerical solutions obtained by different algorithms. Let's consider a PDE system in the operator form $A\widetilde{u} = f$. The discrete operator $A_h$ approximates this system on some grid: $A_h u_h = f_h$. The numerical solution $u_h$ is a grid function (vector $u_h \in R^M$), $\widetilde{u}_h \in R^M$ is the projection of a true solution onto the computational grid, $\Delta\widetilde{u}_h = u_h - \widetilde{u}_h$ is the true approximation error, $\Delta u_h$ is some estimate of this error.

The widely used *a priori* error estimates have the form $\|\Delta\widetilde{u}_h\| \le Ch^l$ ($h$ is the step of discretization, $l$ is the order of approximation) and contain an unknown constant $C$. These estimates are defined for the wide class of solutions and, usually, are not computable for the particular solution. By this reason the *a priori* error estimates have a limited impact on the practical applications and are used mainly at the design and analysis of the numerical algorithms.

The computable measure of the discretization error, which can be used in applications, is commonly denoted as *a posteriori* error estimate $\|\Delta\widetilde{u}_h\| \le E(u_h)$ and is determined by an error indicator $E(u_h)$ that depends on the particular numerical solution $u_h$ and has no unknown constants. In the domain of the finite-element analysis the efficient technique for the *a posteriori* error analysis has the long history (starting from papers [5,6]) and has achieved the maturity at present [7, 8]. It may be successfully applied for the smooth enough solutions (of the elliptic and parabolic equations). In contrast, for the equations of hyperbolic or mixed type (CFD problems at supersonic flow modes, for example) the progress in the *a posteriori* error analysis is limited due to the discontinuities (shock waves, contact surfaces) that may occur and migrate in the flowfield.

Different approaches to the error estimations in CFD are applied at present, which are based on the defect correction [9, 10] or the Richardson extrapolation (RE), for example. The latter is recommended by the modern standards [3,4]. However, the implementation of all these methods reveals significant troubles. In the instance of the defect correction, the precision is limited by the uncertainties occurring at the truncation error estimation due to neglecting the high order terms. The practical applicability of RE is restricted by the variation of the convergence order over the flowfield [11,12]. To settle this issue, the generalized Richardson extrapolation (GRE) is used, which provides an estimate of the error field by the computation of the spatially distributed error order. Unfortunately, GRE may be implemented at the cost of the great computational burden (it requires four or a greater number of the consequent mesh refinements, see [13,14]) and demonstrated the high level of oscillations. There exist computationally cheap approaches for the approximation error norm estimation ([15,16,17,18]), unfortunately, not providing the point-wise information on the error. So, the need exists currently for the computationally inexpensive estimators of the point-wise approximation error in CFD domain for flows with discontinuities.

In this paper we consider the robust and computationally cheap postprocessor-based nonintrusive method ([19]) for the *a posteriori* estimation of the approximation error that uses the set of numerical solutions obtained by the methods of diverse structure on the same grid. The approximation error is computed using the point-wise difference of numerical solutions and the Inverse Problem, formulated in the variational statement with the Tikhonov zero order regularization [20,21]. The numerical tests for two dimensional compressible flows, governed by the Euler equations, are performed to illustrate the features and potentialities of the considered approach. The Edney-I and Edney-VI flow patterns [22] are used due to availability of the analytical solutions. These flow structures are composed by the shock waves and contact discontinuities that provide the most severe conditions for the error estimation. The sets of numerical solutions (three, five, thirteen) computed by distinct numerical algorithms are used for the approximation error estimation. The results are compared with the "true" error obtained by comparing the analytic and numerical solutions.

The paper is organized as follows. In Section 2 we consider the system of linear equations that relates the approximation errors and the differences of the independent numerical solutions. Section 3 recasts the previous Section statement in the form of the variational Inverse Problem with the zero order Tikhonov regularization. Section 4 concerns the quality of considered postprocessor from the viewpoint of the error norm estimation in terms of the effectivity index of the error estimator (in accordance with [7]) and the relative accuracy of the error estimation. The dependence of the computed error on the regularization parameter magnitude is considered in Section 5. The test problems (governing equations and flow patterns) are presented in the Section 6. Section 7 provides the short description of the numerical algorithms used to simulate the above mentioned flows. The results of the numerical experiments for the approximation error estimation are presented in Section 8. The discussion is presented in Section 9. The conclusions are provided in Section 10.

## 2. The approximation errors and differences of numerical solutions

Let's consider a set of grid functions $u_m^{(i)}$ obtained at the discretization and the numerical solution for some system of partial differential equations. We assume some vectorization of multidimensional flowfield data obtained in numerical tests by index $m = 1...M$. Herein, $M$ is the number of the grid points multiplied by the number of flow variables. The algorithms, used for the flowfield simulation, are numbered by $i = 1...n$. The projection of an exact (unknown) solution on the grid is denoted by $\tilde{u}_{h,m}$, the approximation error (also unknown) for $i - th$ solution is denoted by $\Delta u_m^{(i)}$, $u_m^{(i)} = \tilde{u}_{h,m} + \Delta u_m^{(i)}$. This set of solutions is computed on the same grid by $n$ algorithms of the different inner structure (in several cases, by schemes of different approximation order). The differences of two numerical solutions ($i-th$ and $j-th$) may be recast in the following point-wise form:

$$d_{ij,m} = u_m^{(i)} - u_m^{(j)} = \widetilde{u}_{h,m} + \Delta u_m^{(i)} - \widetilde{u}_{h,m} - \Delta u_m^{(j)} = \Delta u_m^{(i)} - \Delta u_m^{(j)}. \tag{1}$$

One may see from the Expression (1) that the computable differences of numerical solutions are equal to the differences of the approximation errors and, hence, contain some information regarding the unknown errors $\Delta u_m^{(i)}$. We treat this information in accordance with the approach described by [19] in order to determine the approximation error $\Delta u_m^{(i)}$. On the set of $n$ numerical solutions one may obtain $N = n \cdot (n-1)/2$ relations casting the following linear system of equations:

$$D_{ij} \Delta u_m^{(j)} = f_{i,m}. \tag{2}$$

Herein $f_{i,m}$ is a vectorized form of the differences $d_{ij,m}$, $D_{ij}$ is the rectangular $N \times n$ matrix, $j = 1...n; i = 1...N$. The summation over the repeating indexes $i, j, k$ is implied starting from this expression. Formally, this system may be resolved for $n$ that is equal (or greater) three. For the simplest case ($n = 3$, $N = 3$) we use the following linear system:

$$\begin{pmatrix} 1 & -1 & 0 \\ 1 & 0 & -1 \\ 0 & 1 & -1 \end{pmatrix} \begin{pmatrix} \Delta u_m^{(1)} \\ \Delta u_m^{(2)} \\ \Delta u_m^{(3)} \end{pmatrix} = \begin{pmatrix} f_{1,m} \\ f_{2,m} \\ f_{3,m} \end{pmatrix}, \tag{3}$$

where

$$\begin{pmatrix} f_{1,m} \\ f_{2,m} \\ f_{3,m} \end{pmatrix} = \begin{pmatrix} d_{12,m} \\ d_{13,m} \\ d_{23,m} \end{pmatrix} = \begin{pmatrix} \Delta u_m^{(1)} - \Delta u_m^{(2)} \\ \Delta u_m^{(1)} - \Delta u_m^{(3)} \\ \Delta u_m^{(2)} - \Delta u_m^{(3)} \end{pmatrix} = \begin{pmatrix} u_m^{(1)} - u_m^{(2)} \\ u_m^{(1)} - u_m^{(3)} \\ u_m^{(2)} - u_m^{(3)} \end{pmatrix}. \tag{4}$$

In the numerical tests we apply the Equation (2) also for five and thirteen variables in similar forms, omitted, herein, for brevity.

Unfortunately, the straightforward solution $\Delta u_m^{(j)} = (D_{ij})^{-1} f_{i,m}$ using the matrix inversion can't be obtained, since the determinant of the matrix $D_{ij}$ should be equal to zero. This degeneracy is caused by the invariance of the solution of system (2) under the shift transformation $\Delta u_m^{(j)} = \Delta \widetilde{u}_m^{(j)} + b$ ($\Delta \widetilde{u}_m^{(j)}$ is the true error) for any $b \in (-\infty, \infty)$ due to the usage of the difference of solutions as the input data. Hence, the solution of the considered problem is nonunique. In accordance with [20], a problem is called well-posed in the Hadamard sense if its solution exists, is unique, and continuously depends on the input data. Otherwise, the problem is called ill-posed. Most Inverse Problems (IP) are ill-posed and by this reason suffer from the generic instability that causes the need for some regularization. The regularized method for the numerical solution of Eq. (2) that was used in the paper is presented in the following Section.

### *3. Inverse Problem statement*

In order to determine $\Delta u_m^{(j)}$ we pose the Inverse Problem in the variational statement [21] with the zero order Tikhonov regularization term:

$$\varepsilon_m(\Delta \bar{u}) = 1/2(D_{ij}\Delta u_m^{(j)} - f_{i,m}) \cdot (D_{ik}\Delta u_m^{(k)} - f_{i,m}) + \alpha/2(\Delta u_m^{(j)} E_{jk} \Delta u_m^{(k)}), \; i = 1...N, j, k = 1...n. \quad (5)$$

Herein, $\alpha$ is the regularization parameter, $E_{jk}$ is the unite matrix. The problem is stated in the point-wise way, separately for each $m - th$ element of vectorized solution (no summation on the repeating $m$).

For the search of $\Delta u_m^{(j)}$, which minimize the functional (5), we use the gradient based method:

$$\Delta u_m^{(j),k+1} = \Delta u_m^{(j),k} - \tau \nabla \varepsilon_m, \quad (6)$$

where $k$ is the iteration number, $\tau$ is a step length parameter.

Since the problem at hand is ill-posed (underdetermined), a regularization should be used to obtain the stable solution. We apply the zero order Tikhonov regularization by the following reason. It is natural to search for solutions of the minimum shift (unavoidable) error $|b|$ (ideally, $|b| \to 0$). For this purpose, we consider the search for the minimal $L_2$ norm of $\Delta u_m^{(j)}$ (normal solution, [20]) denoted as $\delta(b_m)$ in order to stress the dependence on the shift value:

$$\min_{b_m}(\delta(b_m)) = \min_{b_m} \sum_j^n \; (\Delta u_m^{(j)})^2 / 2 = \min_{b_m} \sum_j^n \; (\Delta \tilde{u}_m^{(j)} + b_m)^2 / 2. \quad (7)$$

This relation imposes some restrictions on the magnitude of $|b_m|$. Due to the expression $\Delta\delta(b_m) = \sum_j^n \; (\Delta \tilde{u}_m^{(j)} + b_m)\Delta b_m$ the minimum of (7) over $b_m$ (the norm depends only on $b_m$ at constant true errors) occurs at

$$b_m = -\frac{1}{n}\sum_j^n \; \Delta \tilde{u}_m^{(j)} = -\Delta \bar{u}_m. \quad (8)$$

From this standpoint, the expression (7) describes the minimum of the true error dispersion over the set of solutions (the deviation of the true error from the mean one)

$$\Delta u_m^{(j)} = \Delta \tilde{u}_m^{(j)} + b_m = \Delta \tilde{u}_m^{(j)} - \Delta \bar{u}_m. \quad (9)$$

So, for the problem at hand, the minimum of (7) (which defines the normal solution) corresponds to the shifted true solution. The shift error (8) is equal to the mean true error (with the opposite sign), and, consequently, the accuracy of $\Delta u_m^{(j)}$ estimation is restricted by the unknown mean error value (so, $\Delta u_m^{(j)}$ contains some unremovable error). However, the good point is that the assumption of $\delta(b_m)$ minimality ensures the boundedness of the unremovable error $b_m$. Furthermore, the magnitude of the mean true error $b_m$ should not to be too great. It depends on the errors magnitude and the correlation between them, so it may decay at the enhancement of the ensemble of solutions. The convergence of the regularized solution, obtained by the minimization of (5), to the normal solution is analyzed in [20].

## *4. The quality of the approximation error estimation*

In accordance with [7], the quality of *a posteriori* error estimation may be illustrated by the *effectivity index* of the estimator that is equal to the ratio of the estimated error norm to the true error norm (usually, the global norm is used, which is computed over all grid nodes):

$$I_{eff}^{(j)} = \frac{\left\| \Delta u^{(j)} \right\|_{L_2}}{\left\| \Delta \widetilde{u}^{(j)} \right\|_{L_2}}. \tag{10}$$

One may treat the norms of the true error and estimation error as the radii of hyperspheres $r_{exact} = \left\| \Delta \widetilde{u}^{(j)} \right\|_{L_2}$ and $r_{est} = \left\| \Delta u^{(j)} \right\|_{L_2}$ in the space of discrete solutions. Then the relation $I_{eff}^{(j)} \geq 1$ means that the hypersphere, containing the projection of true error onto grid, belongs to the hypersphere defined by the estimator. In other words, the projection of true solution belongs to some hypersphere with the centre at the approximate solution. So, in order to provide the reliable estimation, the effectivity index should be greater the unit. On the other hand, the estimation should be not too pessimistic, so the value of the effectivity index should be not too great. For the finite elements applications in the domain of elliptic equations (as usual, highly regular), the acceptable range of the effectivity index, according [7], is $1 \leq I_{eff}^{(j)} \leq 3$. It should be mentioned that the upper bound is not strictly defined and may be problem dependent. The solutions, considered herein, contain discontinuities (shear lines, shock waves), so the acceptable range of the effectivity index may be greater and corresponding bounds should be established by additional analysis (from data on the acceptable errors of the valuable functionals used in practical applications).

With account of the expression (9), the effectivity index for $j - th$ solution may be presented as:

$$I_{eff}^{(j)} = \frac{\left\| \Delta \widetilde{u}^{(j)} - \Delta \bar{u} \right\|_{L_2}}{\left\| \Delta \widetilde{u}^{(j)} \right\|_{L_2}}. \tag{11}$$

In dependence on the relation of $\Delta \widetilde{u}^{(j)}$ and $\Delta \bar{u}$, one may obtain the effectivity index value from the range $1 - \left\| \Delta \bar{u} \right\|_{L_2} / \left\| \Delta \widetilde{u}^{(j)} \right\|_{L_2} \leq I_{eff}^{(j)} \leq 1 + \left\| \Delta \bar{u} \right\|_{L_2} / \left\| \Delta \widetilde{u}^{(j)} \right\|_{L_2}$. So, an underestimation of the error ($I_{eff}^{(j)} \leq 1$) is feasible for some numerical solutions belonging to the analyzed set. If $\Delta \widetilde{u}^{(j)} \approx \Delta \bar{u}$ the effectivity index is close to zero and the error estimation fails. Thus, some part of the error estimates may have a low effectivity index. The expansion of the solutions set and the analysis of the distances between solutions (the selection of most remote solutions) may improve the value of the effectivity index. The best estimation $I_{eff}^{(j)} \to 1$ is possible, if $b = -\Delta \bar{u} \to 0$. This convergence may be violated for the numerical methods used in CFD applications due to the correlation of error in vicinity of discontinuities that is caused by the artificial viscosity or limiters. In finite element applications $I_{eff}^{(j)} \to 1$ as the grid step diminishes ($h \to 0$) [7]. Unfortunately, there are no reasons for such asymptotic exactness in the considered approach, since both the difference of solutions and the true error are governed by the terms of the same order of convergence.

The effectivity index is not intuitively lucid, since the negative and positive deviations from the optimal value have the quite different meaning. By this reason, the relative accuracy of the error estimation

$$I_{rel}^{(j)} = \left\| \Delta u^{(j)} - \Delta \widetilde{u}^{(j)} \right\|_{L_2} / \left\| \Delta \widetilde{u}^{(j)} \right\|_{L_2} \tag{12}$$

may be used as another quality indicator for the error estimate. This indicator is more transparent from the intuition viewpoint, since the dependence of the estimate error on the index value is monotonous ($I_{rel}^{(j)} \to 0$ for the precise solution). The values of indicators (10) and (12) are related as $I_{rel}^{(j)} \leq 1 + I_{eff}^{(j)}$ due to the triangle inequality. The relation of these indicators magnitudes reflects a trade between reliability (Eq. (10)) and the accuracy (Eq. (12)) of the estimate.

## 5. The impact of the regularization parameter magnitude on the accuracy of the error estimation

The solution $\Delta u_m^{(j)}(\alpha)$, providing the minimum of the functional (5), depends on the regularization parameter $\alpha$. Due to the ill-posedness, $\Delta u_m^{(j)}(\alpha)$ is not bounded at $\alpha = 0$, hence it is not acceptable at $\alpha = 0$. One may observe the trivial solution $\Delta u_m^{(j)}(\alpha) \to 0$ at $\alpha \to \infty$ that is not acceptable also. According to the practice of the IP regularization [21], there exists a range of $\alpha$ with the weak dependence of the solution $\Delta u_m^{(j)}(\alpha)$ on $\alpha$. In this range $\Delta u_m^{(j)}(\alpha)$ has the minimum deviation from the exact solution $\Delta \tilde{u}_m^{(j)}$ and is adopted as the solution of Inverse Problem.

Figs. 1 and 2 illustrate the influence of the regularization coefficient magnitude on the error estimation for the set of three solutions $\Delta u^{(j)}$, $j = 1,2,3$ (herein, the Eqs. (3-6) are used in the scalar (single point) case). The regularization coefficient runs the values from the interval $\alpha = (10^{-10}, 1)$, true errors $\Delta \tilde{u}^{(j)}$ are equal to (1,-2,3) ($\Delta \bar{u} = 2/3$). The data are provided in the logarithm scale over $\alpha$. Fig. 1 presents the dependence on $\alpha$ for the functional $\varepsilon(\alpha)$ (Expression (5)) marked as "eps", mean sum of errors $\sum |\Delta u^{(i)}(\alpha)|/3$ ("error") and averaged effectivity index (marked as "index"):

$$I_{eff} = \frac{(\sum_{i=1}^{3} (\Delta u^{(i)}(\alpha))^2)^{1/2}}{(\sum_{i=1}^{3} (\Delta \tilde{u}^{(i)})^2)^{1/2}} \,. \tag{13}$$

One may see from the Fig. 1 that the solution diverges at $\alpha \leq 10^{-7}$. The quality of results also deteriorates for $\alpha \geq 1$. This behavior is common for the ill-posed problems. The solution weakly depends on the regularization coefficient in the range $\alpha \in (10^{-6}, 10^{-1})$. Thus, a solution from this range may be accepted as the solution of the considered Inverse Problem.

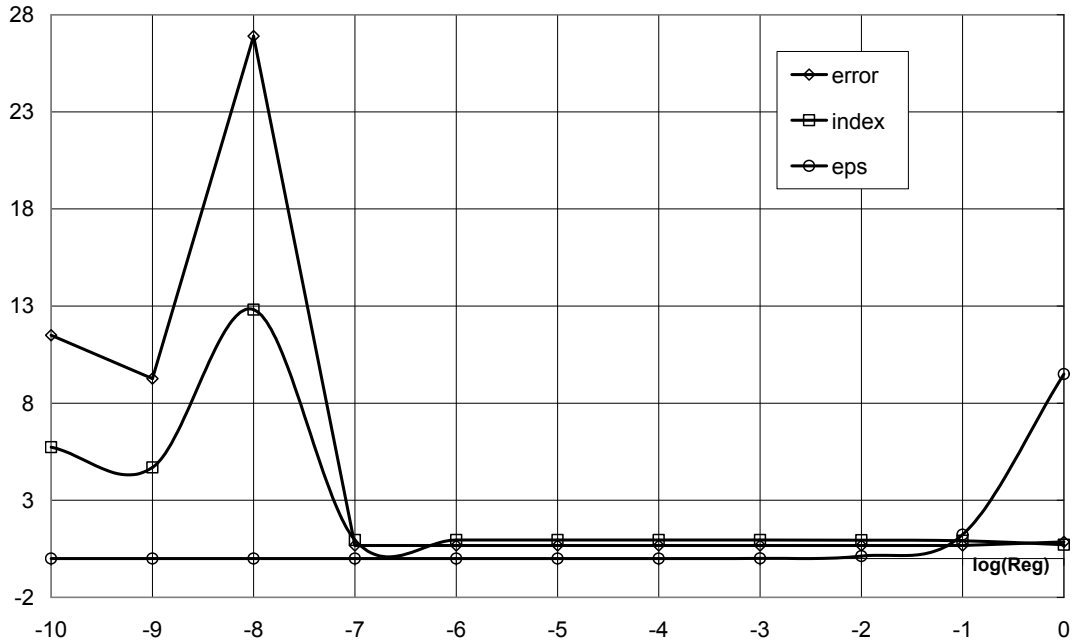

Fig. 1. The dependence of the functional (5) ("eps"), mean error (8) ("error") and the effectivity index (10) ("index") on the regularizing coefficient (in logarithmic scale).

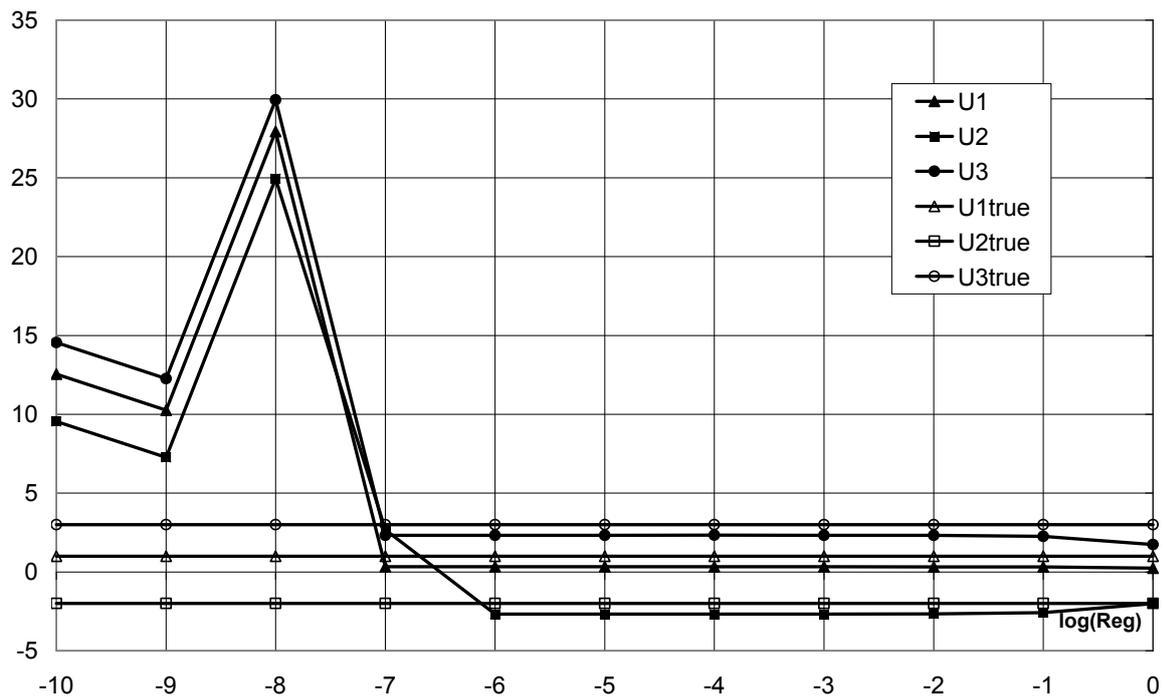

Fig. 2. The dependence of the error estimates on the regularizing coefficient in the comparison with the true errors (in logarithmic scale).

Fig. 2 provides the dependence of the error estimates on the magnitude of the regularizing coefficient $\alpha$ in the comparison with the true errors. The systematic shift $b \approx -0.7$ is observed, which weakly depends on the regularization coefficient and the initial guess. The shift is close to the expected value $b = -\dfrac{1}{3}\sum\limits_{i}^{3}\Delta\widetilde{u}^{(i)} = -2/3$ that corresponds to the average of the true errors $\Delta\overline{u} = 2/3$ (with opposite sign) in this test.

### 6. CFD Test Problems

In order to verify the above analysis and to illustrate the essence of the proposed approach we perform the set of numerical experiments. We consider the approximation error estimation for the tests problems governed by two dimensional compressible Euler equations that follow

$$\frac{\partial\rho}{\partial t} + \frac{\partial(\rho U^k)}{\partial x^k} = 0\,; \tag{14}$$

$$\frac{\partial(\rho U^i)}{\partial t} + \frac{\partial(\rho U^k U^i + P\delta_{ik})}{\partial x^k} = 0\,; \tag{15}$$

$$\frac{\partial(\rho E)}{\partial t} + \frac{\partial(\rho U^k h_0)}{\partial x^k} = 0\,. \tag{16}$$

Here $U^1 = U, U^2 = V$ are the velocity components, $h_0 = (U^2 + V^2)/2 + h$, $h = \dfrac{\gamma}{\gamma - 1}\dfrac{P}{\rho} = \gamma e$, $e = \dfrac{RT}{\gamma - 1}$, $E = \left(e + \dfrac{1}{2}(U^2 + V^2)\right)$ are enthalpies and energies (per unit volume), $P = \rho RT$ is the state equation and $\gamma = C_p/C_v = 1.4$ is the specific heat ratio.

The steady flow patterns, engendered by the interaction of shock waves of I and VI kinds according to Edney classification [22], were used as the test problems due to the availability of analytic solutions (defined by the Rankine-Hugoniot relations on the shock waves). The "true" error is obtained by the subtraction of the numerical solution from the projection of the analytic one on the corresponding grid. It is used for the check of the IP-based estimate of approximation error.

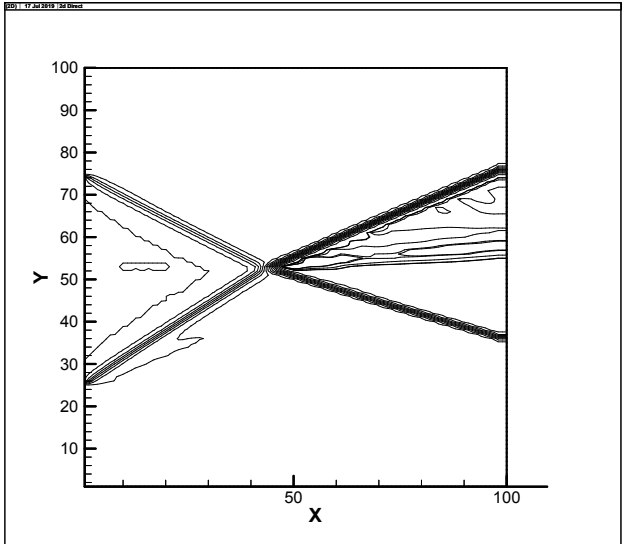

Fig. 3. Density isolines for Edney I flow

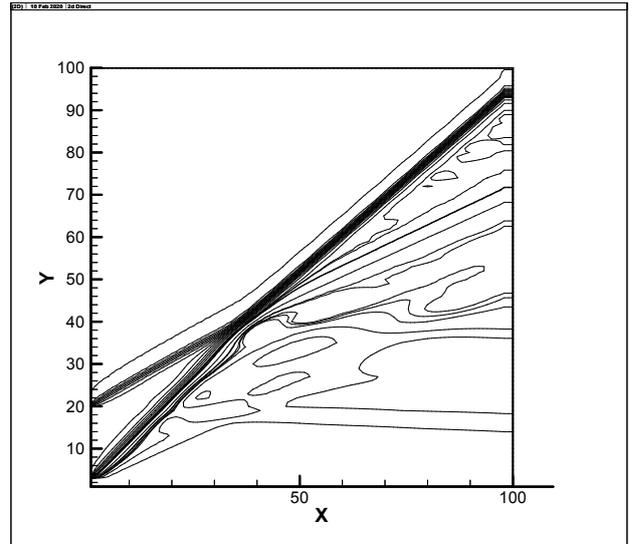

Fig. 4. Density isolines for Edney VI flow

Fig. 3 presents the spatial density distribution for Edney-I flow structure ($M = 4$, $\alpha_1 = 20^o$ and $\alpha_2 = 15^o$). Fig. 4 presents the density distribution for Edney-VI flow structure ($M = 4$, two consequent flow deflection angles $\alpha_1 = 10^o$, $\alpha_2 = 15^o$). Both fields provided in Figs. 3 and 4 are computed using [24].

### *7. The set of used numerical algorithms*

The set of 13 solutions obtained by 10 different numerical algorithms [23-38] and several variants of the artificial viscosity are used in tests in different combinations from the minimal one (three solutions) to maximal (thirteen solutions). The methods and their notations are listed below:

First order algorithm by Courant-Isaacson-Rees [24,25] marked as $CIR$,

Second order MUSCL [26] algorithm based on the approximate Riemann solver by [23] and marked as $AUFS$,

Second order MUSCL [26] algorithm based on the approximate Riemann solver by [27],

Second order algorithm based on the relaxation approach by [28,29],

Second order algorithm based on the MacCormack [30] scheme:

without artificial viscosity,

with second order artificial viscosity with the viscosity coefficient $\mu = 0.01$ and $\mu = 0.002$,

with fourth order artificial viscosity ($\mu = 0.01$),

Second order algorithm based on the "two step Lax-Wendroff" [31,32], with the artificial viscosity of the second order ($\mu = 0.01$),

Third order algorithm based on the modification of Chakravarthy-Osher scheme [33,34],

Third order algorithm WENO [35-37] marked as $W3$,

Fourth order algorithm [38] marked as $ord4$,

Fifth order algorithm WENO [35-37] marked as $W5$.

### *8. The results of numerical tests*

The numerical tests are performed using the uniform grids of $100 \times 100$ and $400 \times 400$ nodes. The results are qualitatively similar and their behavior does not depend on the grid size, so only data for $100 \times 100$ grid are provided by the illustrations.

The Inverse Problem is solved numerically and results are compared with the above discussed "true" error for Edney-I and Edney-VI flow structures. We minimize the functional (5) at every grid point. The gradient of the functional is obtained by the direct numerical differentiation. The regularization coefficient value is $\alpha = 10^{-3}$.

We present the data for error of density $\Delta \rho^{(i)}$ only, since the results for other flow variables are similar.

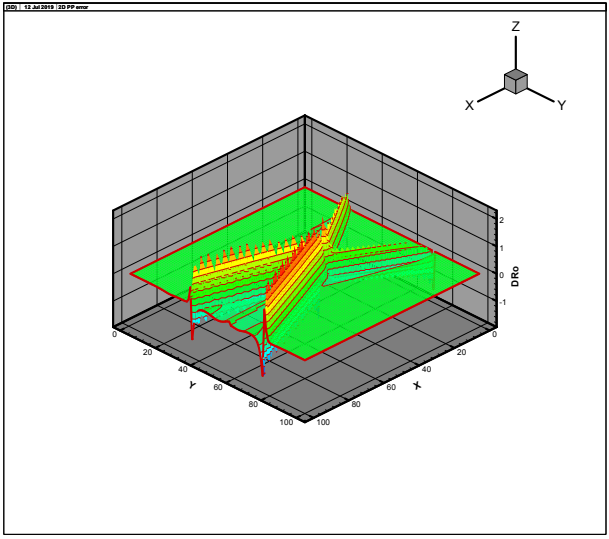

Fig. 5. True error of density for Edney-I flow]

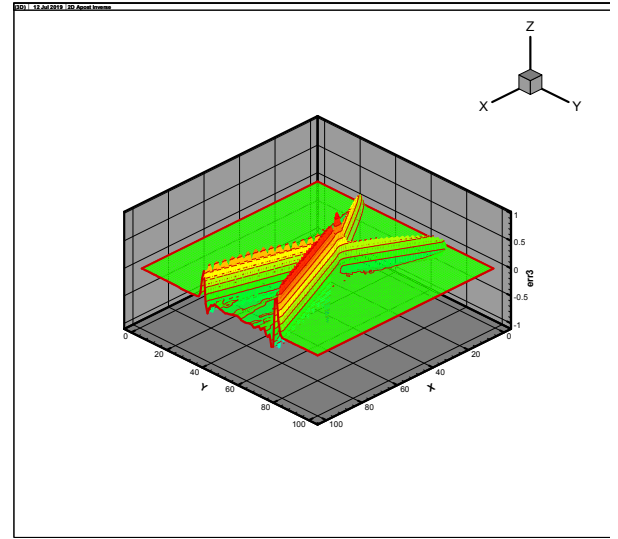

Fig 6. IP estimation of density error for Edney-I flow

Fig. 5 presents the true error for Edney-I flow computed by the first order scheme [24], while Fig. 6 presents the results of the IP based estimation of this error.

Table 1 provides the effectivity index $I_{eff}^{(j)} = \left\| \Delta \rho^{(j)} \right\|_{L_2} / \left\| \Delta \widetilde{\rho} \right\|_{L_2}$ (Eq. (10)) for Edney-I test computed for the different sets of solutions (three solutions are used in two different combinations). Index $j$, which marks the numerical schemes, is replaced by notations, described in the previous Section.

Table 1. *A posteriori* error effectivity index (Eq. (10)) for Edney-I test.

| | $I_{eff}^{CIR}$ | $I_{eff}^{AUFS}$ | $I_{eff}^{ord3}$ | $I_{eff}^{ord4}$ | $I_{eff}^{W5}$ |
|---|---|---|---|---|---|
| three solutions | 0.48 | 0.28 | - | 0.56 | - |
| other three solutions | 0.41 | 0.31 | 0.35 | - | - |
| Five solutions | 1.13 | 1.71 | 1.74 | 1.66 | 1.1 |
| 13 solutions | 0.58 | 1.09 | 0.55 | 0.51 | 0.61 |

The relative accuracy $I_{rel}^{(j)} = \left\| \Delta \rho^{(j)} - \Delta \widetilde{\rho} \right\|_{L_2} / \left\| \Delta \widetilde{\rho} \right\|_{L_2}$ (Eq. (12)) is provided in the Table 2. One may see that the error estimations for the first order method [24] are formally superior due to the relatively great magnitude of the error.

Table 2. The relative accuracy of the error estimation (Eq. (12)) for Edney-I test.

| | $I_{rel}^{CIR}$ | $I_{rel}^{AUFS}$ | $I_{rel}^{ord3}$ | $I_{rel}^{ord4}$ | $I_{rel}^{W5}$ |
|---|---|---|---|---|---|
| three solutions | 0.69 | 1.10 | - | 0.99 | - |
| other three solutions | 0.71 | 1.13 | 1.14 | - | - |
| five solutions | 1.29 | 2.03 | 2.06 | 1.84 | 1.21 |
| 13 solutions | 0.67 | 1.05 | 1.07 | 1.84 | 1.12 |

Fig. 7 presents the true error, while Fig. 8 presents the IP-based estimation of this error for the Edney-VI flow pattern computed by the first order scheme [24].

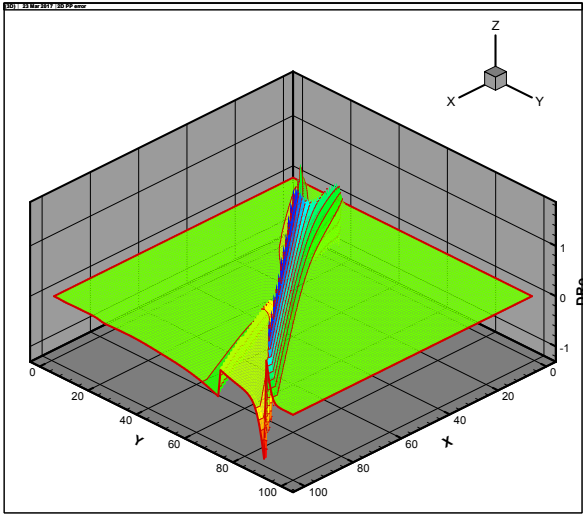

Fig. 7. The true error of density for Edney-VI flow

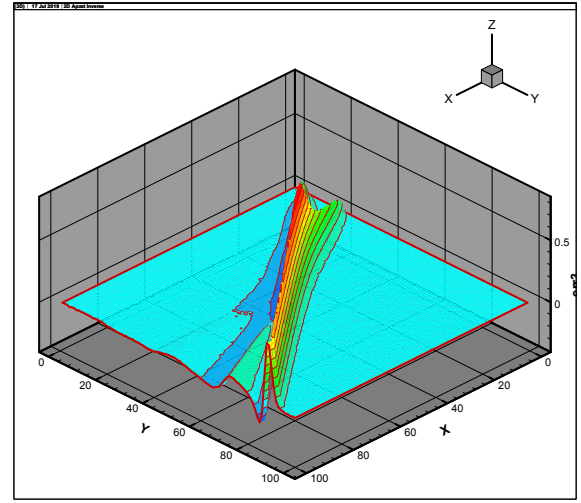

Fig. 8. The IP-based estimation of density error for Edney-VI flow

Table 3 provides the effectivity index (Eq. (10)) for several sets of numerical solutions for Edney-VI test. The relative accuracy (Eq. (12)) is provided by the Table 4 for these solutions.

Table 3. *A posteriori* error effectivity index (Eq. (10)) for Edney-VI test.

|  | $I_{eff}^{CIR}$ | $I_{eff}^{AUFS}$ | $I_{eff}^{ord\,3}$ | $I_{eff}^{ord\,4}$ | $I_{eff}^{W\,5}$ |
|---|---|---|---|---|---|
| three solutions | 2.29 | 3.46 | 3.54 | - | - |
| other three solutions | 2.32 | 3.45 | - | 3.29 | - |
| Five solutions | 0.43 | 0.36 | 0.40 | 0.62 | 0.46 |
| 13 solutions | 0.64 | 0.97 | 0.98 | 0.89 | 0.61 |

Table 4. The relative accuracy of the error estimation (Eq. (12)) for Edney-VI test.

|  | $I_{rel}^{CIR}$ | $I_{rel}^{AUFS}$ | $I_{rel}^{ord\,3}$ | $I_{rel}^{ord\,4}$ | $I_{rel}^{W\,5}$ |
|---|---|---|---|---|---|
| three solutions | 2.38 | 3.62 | 3.70 | - | - |
| other three solutions | 2.38 | 3.61 | - | 3.39 | - |
| five solutions | 0.73 | 1.12 | 1.13 | 1.03 | 0.71 |
| 13 solutions | 0.61 | 0.17 | 0.19 | 0.43 | 0.63 |

For uncorrelated true errors one may expect the improvement of the error estimation quality as the number of solutions is enhanced, since $b_m = -\dfrac{1}{n}\sum_{j}^{n} \Delta \widetilde{u}_m^{(j)}$. The corresponding values of the effectivity index (10) and the relative accuracy (12) are provided in the Tables 1-4 for the sets of 3, 5 and 13 solutions. For the Edney-I test, the reliability (Table 1) and the accuracy (Table 2) of estimates demonstrate no correlated behavior in the dependence on the number of solutions. The accuracy of the estimate increases while the reliability deteriorates as the number of solutions grows. For the Edney-VI test, the improvement of the a posteriori error estimation quality at increase of the number of used solutions is indicated by both indexes. Both the reliability (Table 3) and the accuracy (Table 4) of estimates increase as the number of solutions is enhanced. The uncorrelated and non-monotonic behavior of the reliability (10) and the accuracy (12) indexes at the solutions set expansion may be caused by the non-monotonic convergence to zero of the sum of true errors $b_m$ due to the absence of the ordering of the errors over their norm. The addition of the inaccurate solution may spoil the results obtained at the set of relatively precise solutions. So, in certain situations, the set of three solutions may be optimal for the error estimation from the viewpoint of the computational efforts and the accuracy of results.

Fig. 3 and Fig. 4 provide the comparison of the true error with the results of the IP based error estimation for Edney-I flow. The corresponding analysis for the Edney-VI flow pattern may be performed using Figs. 5 and 6. The errors in the vicinity of shocks have the wave-like shape with the positive and negative half-waves that is quite natural for monotonous smoothing of the stepwise discontinuities. From the qualitative viewpoint Figs. 3-6 show that the estimates of the error demonstrate the similar spatial structure if compare with the true error. The quantitative comparison of the estimates and true errors may me obtained using Figs. 9 and 10, which present the pieces of vectorized grid function of density error obtained by the Inverse Problem in comparison with the true error for Edney-I test. The index along abscissa axis is related with the indexes along $X$ $(k_x)$ and $Y$ $(m_y)$ as $i = N_x(k_x - 1) + m_y$.

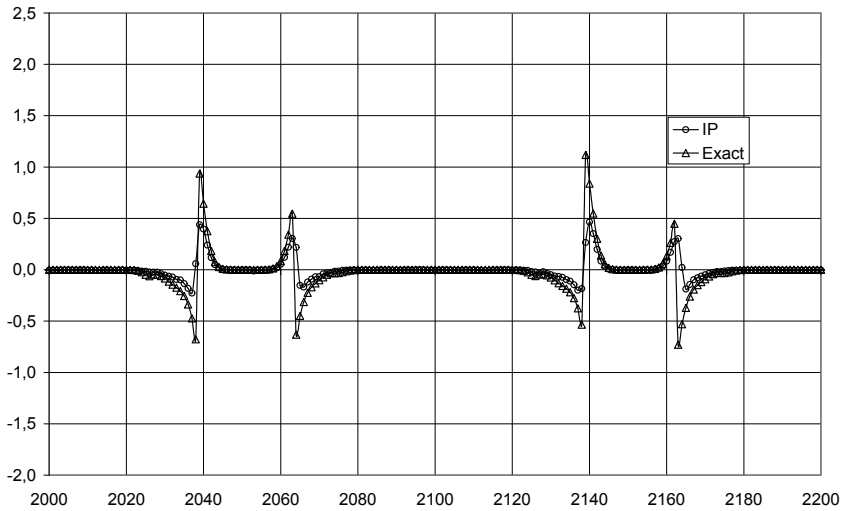

Fig. 9. The comparison of the density error, estimated by the Inverse Problem, with the true error in zone before shocks crossing (Edney-I)

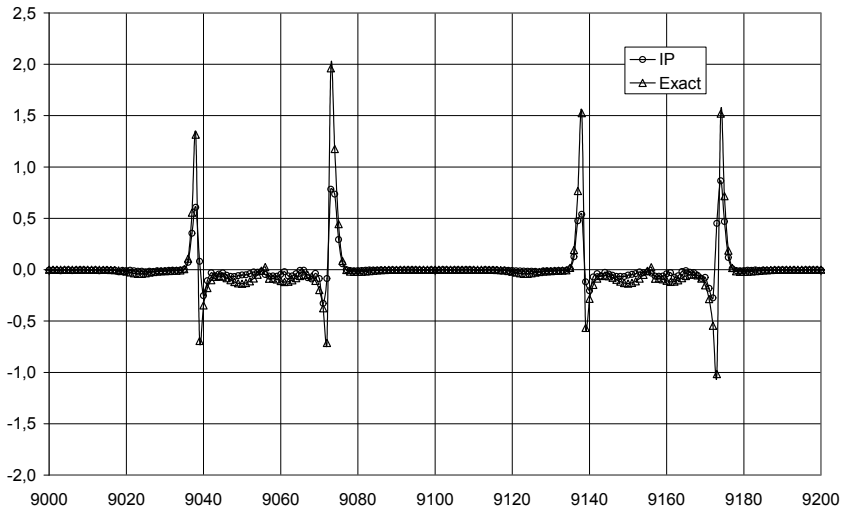

Fig. 10. The comparison of the density error, estimated by the Inverse Problem, with the true error in zone past shocks crossing (Edney-I)

The periodical jump of the flow parameters corresponds to the transition through the shock waves, herein two transitions are presented. Fig. 9 corresponds the domain before shock crossing, Fig. 10 illustrates the flow structure past shocks crossing.

### 9. Discussion

The above considered Inverse Problem based method is less accurate if compared with the Richardson extrapolation due to the presence of the irremovable error, proportional to the mean error over the set of used solutions. Nevertheless, the generalized Richardson extrapolation, which may be correctly used for the shocked flows, suffers from the oscillations due to the variable convergence order that may cause severe point-wise errors [13,14]. From this standpoint, the IP based method is more robust. Additionally, the considered postprocessor is much more computationally inexpensive, since it applies the computations without a mesh refinement. Recall that the calculations are performed on the same grid for all considered numerical methods.

The results by [19], obtained for the supersonic flows over the cone, demonstrates the successful error estimation using three or five solutions, computed by the different algorithms of the same (second) order and certain increase of the quality of results at the expansion of the set. The present paper provides the successful error estimation for the flows engendered by the interaction of the shock waves, if at least three numerical solutions, obtained by methods of different approximation orders, are available. In general, the extension of the number of solutions enables to increase the reliability and accuracy of the results, unfortunately, in the nonmonotonous manner.

Despite the provided tests are selected from the CFD domain, the considered method may be used for the numerical solution of any PDE system without the loss of generality, since no assumptions, which are specific for CFD, are used.

### 10. Conclusions

The set of numerical solutions obtained using different algorithms provides the information for the estimation of the approximation error, which may be computed using the Inverse Problem in the variational statement with the zero order Tikhonov regularization.

Since the considered problem is underdetermined, the results contain the unremovable error, which is equal to the true error, averaged over the ensemble of solutions.

The numerical tests demonstrate the feasibility for the estimation of the pointwise approximation error with the acceptable accuracy for the compressible inviscid flows, engendered by the shock waves interference.